\documentclass[journal]{IEEEtran}
\usepackage{makeidx}
\usepackage{graphicx}
\usepackage{mathptmx}
\usepackage{amsmath}
\usepackage{latexsym}
\usepackage{epsfig}
\usepackage{color}
\usepackage{algorithm}
\usepackage{algorithmicx}
\usepackage{enumerate}
\usepackage[colorlinks,urlcolor=blue,citecolor=blue,linkcolor=blue]{hyperref}

\begin{document}
%
\title{Contingency-Risk Informed  Power System Design}
%
%
%
\author{Richard Li-Yang Chen~\IEEEmembership{Member,~IEEE}, Amy Cohn, Neng Fan~\IEEEmembership{Member,~IEEE}, and Ali Pinar~\IEEEmembership{Senior Member,~IEEE}
\thanks{R. Chen is with Quantitative Modeling and Analysis, Sandia National Laboratories, Livermore, California 94551, USA. email: rlchen@sandia.gov.}
\thanks{A. Cohn is with Industrial and Operations Engineering, University of Michigan, Ann Arbor, Michigan 48109, USA. email: amycohn@umich.edu.}
\thanks{N. Fan is with Systems and Industrial Engineering Department, University of Arizona, Tucson, AZ 85721, USA. email: nfan@email.arizona.edu.}
\thanks{A. Pinar  is with Quantitative Modeling and Analysis, Sandia National Laboratories, Livermore, California 94551, USA. email:apinar@sandia.gov }}

%
%

\markboth{IEEE Transactions on Power Systems}%
{Contingency Risk Aware Power System Design}

\maketitle

\begin{abstract}
We consider the problem of designing (or augmenting) an electric power system at a minimum cost  such that it satisfies the $N$-$k$-$\boldsymbol \varepsilon$ survivability criterion. This survivability criterion is a generalization of the well-known $N$-$k$  criterion, and it requires that at least $(1- \boldsymbol \varepsilon_j)$  fraction of the total demand to be met after failures of  up to   $j$  components, for $j=1,\cdots,k$.  The network design problem adds another level of complexity to the notoriously hard  contingency analysis problem, since the contingency analysis is only one of the requirements for the design optimization problem.
We present a mixed-integer programming formulation of this problem  that takes into account both transmission and generation expansion. We propose an algorithm that can avoid combinatorial explosion in the number of contingencies, by seeking vulnerabilities in intermediary solutions and constraining the design space accordingly.
Our approach is built on our  ability to identify such system vulnerabilities quickly.   Our empirical studies  on modified instances from the IEEE 30-bus and IEEE 57-bus systems  show the effectiveness of our methods. We were able to solve the transmission and generation expansion problems for $k=4$ under 2 minutes, while other approaches failed to  provide a solution at the end of 2 hours.

\end{abstract}


%
\IEEEpeerreviewmaketitle

\section{Introduction}\label{sec1}
\IEEEPARstart{C}{ontingency} analysis of power systems is a problem of increasing importance due to ever increasing demand for electricity. This increase in electricity demand outpaces  the increase in  system capacity growth, which  leaves less  room for redundancy in the system.  As a consequence, equipment failures are now more likely to lead to blackouts.  At the same time,  society's increased reliance on electricity breeds  far and wide ramifications for a blackout. While  minimizing the likelihood of any blackout is a laudable goal, this will inevitably come at an economic cost.   A partial remedy for increasing costs is  to take into account both the likelihood of an event  and  severity of its consequences in the planning phase. This concept   is adopted  in Transmission Planning Standard (TPL-001-1, \cite{NERC}), defined by the North American Electric Reliability Corporation (NERC).   According to this standard, if  only a single element is lost ($N-1$ contingency), the system must be stable, without any loss-of-load. In the case of $k$ simultaneous  failures ($N-k$ contingency), the system still has to restore stability, but a limited loss-of-load is allowed.

While contingency analysis, as a concept, is simple,  algorithmically it is far from trivial.  Despite its difficulty,  the critical importance of the problem has drawn a lot of interest, and significant progress has been made over the last decade. In particular,  optimization methods have been proposed to replace enumerative approaches that rely on  verifying feasibility  on each state of a given list of contingencies.
For example, optimization methods are used to find small groups of lines whose failure can cause a severe blackout or a large loss-of-load \cite{Donde08, lesieutre08, Pinar2010,Bienstock2010,Arroyo2010}.  $N$-$k$ contingencies were also studied in optimal power flow models \cite{SalWB04,SalWB09,Fan2011} and unit commitment problems \cite{Street2011} using a single bus model.
The methods used in \cite{Arroyo2010,Fan2011,Street2011} are all based on  a bilevel programming approach, which is the main method used for network inhibition/interdiction problems.
Probability analysis \cite{CheM2005}, limitations on generation and line capacity \cite{CarL2002}, swarm optimization \cite{LiC2009}, and methods based on topological characteristics of power grids \cite{HinC2010} have been used for vulnerability analysis.
Other approaches for contingency analysis are reviewed in \cite{MorG2001,CheJ2009}.

The great strides in contingency analysis over the last decade have established the basis for higher objectives. In this paper, we consider the transmission and generation expansion planning (TGEP) problem with contingency constraints.  Our goal is to design (or augment) a power system at a minimum cost  to satisfy contingency constraints.  We  formalize the contingency requirements with  the $N$--$k$--$\boldsymbol \varepsilon$ criterion, where  $\boldsymbol \varepsilon$ is a parameter \emph{vector} that specifies allowable loss-of-load, for each contingency size, as a fraction of total system demand.  More specifically, this criterion requires that for any contingency of size $j=0,1,\cdots,k$,  at least $(1-\varepsilon_j)$ fraction of the total demand must be satisfied. Following NERC's TPL-001-1 standards, for the no-contingency state and contingencies of size one $(j=1)$, no loss-of-load is allowed (i.e. $\varepsilon_0=\varepsilon_1=0$); for multiple failure contingencies $(j\ge 2)$, a small fraction of total load demand can be shed (i.e. $0 < \varepsilon_2 \le \cdots \le \varepsilon_k < 1$).

To understand the complexity of the TGEP, one should observe that contingency analysis, a difficult combinatorial optimization problem by itself, is only one of the prerequisite steps in solving TGEP; we must also address network design issues.  Our approach is built on  our ability to solve the contingency analysis problem efficiently so that it can be embedded in a broader framework.
We begin by formulating a mixed-integer linear program (MILP) to model TGEP that explicitly includes multiple states representing \emph{each} of the possible contingency states and, for each of these contingency states, contains the corresponding generation and flow variables to ensure that at least $(1-\varepsilon_j)$ fraction of the demand can be met. Such an enumerative approach, however, becomes computationally intractable even for  moderate values of system size $N$ and maximum contingency-size budget $k$, since the number of states grows with $N \choose k$.
Furthermore, it should be noted that contingency analysis for TGEP requires looking at higher values of $k$, compared with contingency analysis for  shorter term problems.   The reason  is that in long term problems, planned outages, such as  system components going under maintenance,  contribute to $k$, whereas in shorter-term problems they do not.  For instance,  for day-ahead  planning,  we will know the unavailable system components in advance; remove them from the system; and then we will do contingency analysis for unexpected failures in the remaining system. However, for long term problems such as TGEP, we want to design a system such that a feasible direct current optimal power flow (DCOPF) exists for any combination of planned and unplanned outages.  This forces us to look at higher values of $k$, which  makes enumerative approaches prohibitively expensive, even when state of the art high performance computing platforms are available.

To overcome this challenge, we propose two cutting plane algorithms, one based on a direct application of Benders decomposition method to check the load satisfaction of each contingency state \emph{explicitly} and another based on an online contingency state generation (OCS) algorithm, which solves bilevel separation problems to determine the worst-case loss-of-load for each contingency size $j=0,1,\cdots,k$. The OCS  algorithm  \emph{implicitly} identifies worst-case contingencies without explicitly evaluating all contingency states, and generates additional constraints, corresponding to Benders feasibility constraints, to exclude solutions that are infeasible under the identified worst-case contingencies.  This approach  avoids the prohibitive cost of generating a constraint for every combination of failures, by only identifying relevant combinations as needed.  As the computational experiments show,  this reduction to only relevant contingencies makes a tremendous impact on the scalability of the algorithm.

 We applied  these approaches to  the IEEE 30-bus and IEEE 57-bus systems.  Computational results show the scalability of the OCS algorithm. We were able to solve the TGEP problem for up to $k=4$ in under 2 minutes, while explicit enumeration approaches, e.g. extensive form and Benders decomposition,  failed  to provide any solutions at the end of 2 hours. We observed the key to our success was our ability to avoid  looking at all individual vulnerabilities, and only a small number of iterations  (vulnerability searches) are enough to find a provably-optimal solution.

Transmission  and generation expansion problems  have been studied for a long time. A survey of earlier work in this area can be found in~\cite{Latorre2003}. Recently, especially after the 2003 Northeast American blackout, security, stability and reliability issues have been becoming another major concerns in modern power systems.  Our earlier results were presented in~\cite{nsw11, Pmaps11}.
%
 \cite{Delgadillo2011} and \cite{Romero2012} addressed the defense protection of power system, and they analyzed the interaction between a power system defender and a terrorist who seeks to disrupt system operations. \cite{Carrion2007}, \cite{Choi2007},  and \cite{Moulin2010} studied the contingency criteria by stochastic programming and integer programming approaches.  \cite{Zhang2011} and \cite{Romero2012} proposed a multilevel mixed integer programming models for transmission expansion. Most of the work in this area however, is restricted to the $N-1$ contingency criteria, which cannot reflect current situation for contingencies with more than one failure. Additionally, the proposed models cannot be directly extended for this new situation, and they only consider the transmission expansion, and contingencies are restricted to transmission elements.

The generation expansion problem has recently been studied in \cite{Jin2011}.
Transmission expansion and generation expansion planning problems have also been studied in a unified  model \cite{Bent2011}.
Integration of  renewable energy resources was taken into account  generation and transmission expansion planning~\cite{Weijde2010,Bent2011,Jin2011}. However, none of these studies  consider the contingency criteria in case of failures of both generating unit and transmission elements.



The rest of this paper is organized as follows. In Section \ref{sec2}, the TGEP problem considering the full set of contingency states is formulated as a mixed integer nonlinear program (MINLP); Section \ref{sec3} presents two methods to solve this large-scale MINLP; in Section \ref{sec4}, results of numerical experiments performed on two IEEE test systems are presented; Section \ref{sec5} concludes the paper.

\section{Models}\label{sec2}
\subsection{Nomenclature}
\noindent {\bf Sets and indices}
\begin{tabbing}
\hspace{6ex}\=\hspace{6ex}\=\hspace{4ex}\=\hspace{4ex}\=\hspace{4ex}\=\hspace{4ex}\=\hspace{4ex}\=\hspace{4ex}\kill
$V$ \>      Set of buses (indexed by $i,j$).\\
$S_j$\>  Set of \emph{all} contingency states with \emph{exactly} $j$ failures\\
\>(indexed by $s$).\\
$S$ \> Set of \emph{all} contingency states with $k$ \emph{or less} failures,\\
\>$S = S_1\cup S_2 \cdots \cup S_k$.\\
$|s|$ \> Number of failed element(s) in contingency $s$.\\
$G$ \>  Set of generating units (indexed by $g$).\\
$G_i$\> Set of generating units at bus $i$.\\
$E$ \>  Set of transmission elements (indexed by $e$).\\
$E_{.i}$ \>  Set of transmission elements oriented into bus $i$.\\
$E_{i.}$ \>  Set of transmission elements oriented out of bus $i$.\\
$i_e,j_e$ \>    Tail/head (bus no.) of transmission element $e=(i_e,j_e)$.\\
\end{tabbing}

\noindent {\bf Parameters}
\begin{tabbing}
\hspace{6ex}\=\hspace{6ex}\=\hspace{4ex}\=\hspace{4ex}\=\hspace{4ex}\=\hspace{4ex}\=\hspace{4ex}\=\hspace{4ex}\kill
$C_e$ \>     Investment cost of transmission element $e$.\\
$C_g$\>      Investment cost of generating unit $g$.\\
$C^p_g$ \>   Marginal production cost of generating unit $g$.\\
$\overline{P}_g$ \>   Maximum capacity of generating unit $g$.\\
$B_e$ \>     Electrical susceptance of transmission element $e$.\\
$F_{e}$ \>   Capacity of transmission element $e$.\\
$D_i$ \>     Electricity load demand at bus $i$.\\
$D$ \>     Total electricity load demand across all bus $i \in V$,\\
\>  $D=\sum_{i \in V}D_i.$\\
$M_e$ \>     Big M value for transmission element $e$.\\
$\sigma$ \>  Weighting factor to make investment cost and \\
\> operating cost comparable.\\
$N$ \> Total number of transmission elements and \\
\>generating units, $N=|G|+|E|$.\\
$k$ \> Maximum contingency size under consideration. \\
\>In any contingency state the number of  failures, i.e. the  \\
\>number of system elements in a contingency, is  \\
\>between  $0$ and $k$.\\
$\varepsilon_j$\>    Fraction of load demand that can be shed given size $j$\\
\> contingencies, for all $j=0,1,\cdots,k$.\\
$\tilde d_g^s$ \>        Binary parameter that is 1 if generating unit $g$ is \\
\> part of the contingency state $s$ and 0 otherwise.\\
$\tilde d_e^s$ \>        Binary parameter that is 1 if transmission element $e$ \\
\>                is part of the contingency state $s$ and 0 otherwise.\\
$\tilde {\textbf d}^s$ \> Vector that  concatenates $\tilde d_g^s$  and $\tilde d_e^s$ variables.\\
\end{tabbing}

\noindent {\bf Decision variables}
\begin{tabbing}
\hspace{6ex}\=\hspace{6ex}\=\hspace{4ex}\=\hspace{4ex}\=\hspace{4ex}\=\hspace{4ex}\=\hspace{4ex}\=\hspace{4ex}\kill
$x_g$ \>  Binary variable that is 1 if generating unit  $g$ is added \\
\>          and 0 otherwise.\\
$x_e$ \>  Binary  variable that is 1 if  transmission element $e$ is \\
\>           added and 0 otherwise.\\
$\textbf x$ \>  Vector that  concatenates $x_g$  and $  x_e$ variables.\\
$f_e$ \>  Power flow on transmission element $e$.\\
$p_g$ \>  Power output of generating unit $g$.\\
$q_i$ \>  Loss-of-load at bus $i$.\\
$\theta_i$ \> Phase angle of bus $i$.\\
$d_g$ \>        Binary variable that is 1 if generating unit $g$   is \\
\>        \emph{selected}     to be    in the contingency and 0 otherwise.\\
$d_e$ \>        Binary variable that is 1 if transmission element $e$ is \\
\>               \emph{selected} to be in the contingency and 0 otherwise.\\
$f^{s}_e$ \>  Power flow  for transmission element $e$ for \\
\> contingency state $s$.\\
$p^{s}_g$ \>  Power output of generating unit $g$ for contingency \\
\>state $s$.\\
$q^{s}_i$ \>  Loss of load at bus $i$ for contingency state $s$.\\
$\theta_i^{s}$ \> Phase angle of bus $i$ for contingency state $s$.
\end{tabbing}

For a contingency state $s\in S$, $\tilde{d}_g^s=1$ and $\tilde{d}_e^s=1$ denote that generating unit $g$ and transmission element $e$ fail in this state, respectively. Conversely, $\tilde{d}_g^s=0$ and $\tilde{d}_e^s=0$ denote that both these two elements are available. Thus, for $s=0$ denoting the no-contingency state,  we have $\tilde{d}_g^0=0$ and $\tilde{d}_e^0=0$ for all $g\in G$ and $e\in E$. 
%

\subsection{Transmission and Generation Expansion  Model}
In this section, we extend the standard TGEP problem to include contingency constraints. For brevity of presentation, we treat all power system elements, including existing ones,  as candidates for addition to the system.  For an existing element, the investment cost $(C_e, C_g)$ is set to $0$ and the corresponding investment decisions $(x_e, x_g)$ fixed at $1$.

Once the network design decisions are made, each element selected for addition becomes available in all contingency states $s\in S_j$ and $j=0,1,\cdots,k$, unless it is part of a given contingency. Consistent with the NERC reliability standards, in the no-contingency state ($N$-$0$) and the single contingency state ($N$-$1$), no loss-of-load is allowed, thus $\varepsilon_0=\varepsilon_1=0$. For larger contingency states, e.g. $j = 2,\cdots, k$, the total loss-of-load is limited by the threshold $\varepsilon_j$.  We assume that $\varepsilon_j \leq \varepsilon_{j+1}$ for all $j=0,1,\cdots, k-1$. The MINLP model for TGEP is formulated as follows,

\begin{subequations}\label{ef}
\begin{align}
\min_{\tilde {\mathbf x}, \mathbf f,\mathbf p, \mathbf q, \boldsymbol \theta}  \ &\sum_{e \in E} C_e x_e + \sum_{g\in G} C_gx_g + \sigma  \sum_{g\in G}  C_g^p p_g^0 \label{ef_obj} \\
\textrm{s.t.} \quad
	&\sum_{g\in  G_i} p_g^{s}  + \sum_{e \in E_{.i}} f_{e}^{s} -\sum_{e \in E_{i.}} f_{e}^{s} + q_i^s = D_i, \quad  \forall i \in V, \label{ef_node_bal} \\
 &\hskip 5.8cm s \in S  \nonumber \\
    &B_{e}\big(\theta_{i_e}^{s}- \theta_{j_e}^{s}\big) x_e (1-\tilde d_e^s) - f_{e}^{s}=0, \quad \forall e \in E, \label{ef_flow_cons}\\
 &\hskip 5.8cm s \in S \nonumber \\
	&-F_e x_e(1-\tilde d_e^s)\leq f_{e}^{s} \leq  F_{e}x_e(1-\tilde d_e^s), \quad \forall e \in E, \label{ef_flow_cap} \\
 &\hskip 5.8cm s \in S \nonumber \\
	&0 \leq p_g^{s} \leq \overline{P}_g x_g(1-\tilde d_g^s), \quad \forall g \in G, s \in S\label{ef_p_bounds}\\
    &\sum_{i \in V}  q_i^s \le \varepsilon_{|s|} D,  \quad \forall s\in S \label{ef_lolc} \\
        &0 \leq q^{s}_i \leq D_i, \quad \forall i \in V, \forall s \in S\label{ef_q_lb}\\
            	&x_g \in \{0,1\}, \quad \forall g \in G \label{ef_x}\\
	&x_e \in \{0,1\}, \quad \forall e \in E \label{ef_x2}
\end{align}
\end{subequations}
In all subsequent formulations, unless otherwise specified, the indices $i,g,e,s$ and $k$ are elements of their corresponding sets, i.e., $i\in V, g\in G, e\in E$, and $s\in S$.

The objective \eqref{ef_obj} is to minimize the total transmission and generation investment cost plus  the weighted operating cost in the no-contingency state $(s=0)$.
Note that our formulation does not take into account the operational costs during a contingency state. There are two reasons for this. First, contingencies have low likelihood and thus the operational costs during such events are negligible. The real financial burden of a contingency shows itself when the contingency leads to a blackout, which brings us to the second reason: the primary goal during a contingency is to keep the system intact as opposed to minimizing operational costs, since the cost of system failure is likely to be significantly  higher than any operational cost.  Constraints \eqref{ef_node_bal} are conservation of flow requirements for each bus and contingency pair. For any transmission element that is operational, Kirchhoff's voltage law must be enforced by \eqref{ef_flow_cons}. Power flow on transmission element $e$ is governed by thermal capacity constraints \eqref{ef_flow_cap}. For each contingency state, the power output of a generating unit must satisfy the upper bound given by \eqref{ef_p_bounds}.  Constraints \eqref{ef_lolc} dictate that the total loss-of-load in the system cannot exceed $\varepsilon_{|s|}$ fraction of the total system load, i.e., at least $(1-\varepsilon_{|s|})D$ of demand must be satisfied. Constraints \eqref{ef_q_lb} restrict the loss-of-load at each bus to be at most the total demand at the bus.

Observe that constraints \eqref{ef_node_bal}-\eqref{ef_q_lb} are specific to a particular contingency state; that is, for a given contingency state $s$, the transmission and generation element(s) in the contingency have zero capacity.  In the no-contingency state $s=0$, all invested transmission elements and generating units are available for the DCOPF problem.

\section{Solution Approaches}\label{sec3}
Replacing constraints \eqref{ef_flow_cons} by
\begin{align}
&B_{e}\big(\theta_{i_e}^{s} - \theta_{j_e}^{s}\big) - f_{e}^{s} + M_e(1-x_e +\tilde d_e^s) \geq  0, \ \forall e, s, \\
&B_{e}\big(\theta_{i_e}^{s} - \theta_{j_e}^{s}\big) - f_{e}^{s} - M_e(1 -x_e + \tilde d_e^s) \leq  0, \ \forall e, s.
\end{align}
where $M_e$ is a sufficiently large constant, transforms formulation \eqref{ef} to a mixed integer linear program (MILP), which we refer to as the \emph{extensive form} (EF).  EF will typically have an extremely large number of variables and constraints because it grows with the number of contingency states, which increases exponentially with $N$ and $k$.  For large power systems and/or  a contingency $k$ greater than one, EF rapidly becomes computationally intractable for increasing system size, $N$, and increasing contingency size, $k$.  In the following sections, we modify this formulation and present cutting plane algorithms for solving the reformulated problem.

\subsection{Benders Decomposition}
We begin by presenting an alternative formulation with only $|G|+|E|$ binary variables but possibly an extremely large number of constraints.  We use linear programming duality to generate valid inequalities for the projection of the natural formulation onto the space of the ${\tilde {\mathbf x}}$ variables.  In essence, we use a variant of Benders Decomposition, in which we only generate valid inequalities corresponding to ``feasibility'' cuts.

For a given contingency state, $s\in S$, capacity expansion vector $\tilde {\mathbf x}$ and contingency state vector $\tilde {\mathbf d}^s$, we solve the following linear program, denoted as the \emph{primal subproblem} PSP$(\tilde {\mathbf x},\tilde {\mathbf d}^s)$, to determine a DCOPF that minimizes total system loss-of-load.

\begin{subequations}\label{psp}
\begin{align}
    \min_{\mathbf f,\mathbf p, \mathbf q, \boldsymbol \theta} \quad & \sum_{i \in V} q_i^s \label{psp_obj} \\
	 \textrm{s.t. }     
	(\alpha_i^s)       \quad & \sum_{g\in  G_i} p_g^s+ \sum_{e \in E_{.i}} f_{e}^{s}-\sum_{e \in E_{i.}} f_{e}^{s}  +  q_i^s = D_i,  \ \forall i\label{psp_flow_consv} \\
	(\hat \beta_e^s)   \quad &-B_{e}\big(\theta_{i_e}^{s} - \theta_{j_e}^{s}\big) + f_{e}^{s}  \leq M_e(1- \tilde x_e + \tilde d_e^s), \ \forall e   \\
	(\check \beta_e^s) \quad &B_{e}\big(\theta_{i_e}^{s} - \theta_{j_e}^{s}\big) - f_{e}^{s} \leq M_e(1- \tilde x_e + \tilde d_e^s), \ \forall e   \\
	(\delta_e^s)       \quad &f_{e}^{s} \leq  F_{e}  \tilde x_e(1-\tilde d_e^s), \ \forall e     \label{psp_flow_ub} \\
	(\eta_e^s)         \quad &-f_{e}^{s} \leq  F_{e}  \tilde x_e(1-\tilde d_e^s), \ \forall e     \label{psp_flow_lb} \\
	(\zeta_g^s)        \quad & 0  \leq p_g^{s} \leq \overline{P}_g  \tilde x_g(1-\tilde d_g^s), \ \forall g  \label{psp_gen_ub} \\
	(\lambda_i^s)      \quad & 0  \leq q_i^{s} \leq D_i,\ \forall i   \label{psp_lol_ub}
\end{align}
\end{subequations}

In this formulation, the objective, \eqref{psp_obj}, is to minimize total  loss-of-load by adjusting the flow, phase angles and power generation, given the prescribed capacity expansion decision $\tilde {\mathbf x}$ and contingency vector $ \tilde{\mathbf d}^s$, corresponding to scenario $s$.  Letting $z(\tilde {\mathbf x}, \tilde{\mathbf d}^s)$ be the optimal objective value of (\ref{psp}), if $z(\tilde {\mathbf x}, \tilde{\mathbf d}^s)>\varepsilon_{|s|}D$, there does not exist a feasible DC power flow satisfying at least $(1-\varepsilon_{|s|})$ of total demand.  Alternatively, if $z(\tilde {\mathbf x}, \tilde{\mathbf d}^s) \le \varepsilon_{|s|} D$, there exists a power flow that can satisfy \emph{at least} $(1-\varepsilon_{|s|})$ of total demand.

Variables in parentheses on the left-hand-side of the constraints in \eqref{psp} denote the corresponding dual variables. In turn, we can formulate the dual of this problem, DSP$(\tilde {\mathbf x},\tilde {\mathbf d}^s)$ as follows,
\begin{align*}
\max_{\boldsymbol \alpha, \hat \beta, \check \beta, \delta, \eta, \zeta, \lambda} \ \sum_{i \in V} D_i (\alpha_i^s + \lambda_i^s) + \sum_{e \in E} M_e(1- x_e + \tilde d_e^s)(\hat \beta_e^s + \check \beta_e^s) \\
+ \sum_{e \in E} F_e  x_e(1-\tilde d_e^s)(\delta_e^s + \eta_e^s) + \sum_{g \in G} \overline{P}_g   x_g(1-\tilde d_g^s)\zeta_g^s,  \nonumber 
\end{align*}
subject to constraints corresponding to primal variables $\mathbf f, \mathbf p, \mathbf q, \theta$. Since PSP$(\tilde {\mathbf x},\tilde {\mathbf d}^s)$ has a finite optimal solution value (in the worst case, all load will be shed), DSP$(\tilde {\mathbf x},\tilde {\mathbf d}^s)$ also has a finite optimal solution value and by strong duality, the optimal solutions coincide. Since DSP$(\tilde {\mathbf x},\tilde {\mathbf d}^s)$  has a finite optimal solution value it also has an optimal extreme point. Thus we can reformulate PSP$(\tilde {\mathbf x},\tilde {\mathbf d}^s)$ as follows,
\begin{align}\label{dsp2}
&\max_{\ell\in L^s}  \sum_{i \in V} D_i (\alpha_i^{\ell} + \lambda_i^{\ell}) + \sum_{e \in E} M_e(1- \tilde x_e + \tilde d_e^s)(\hat \beta_e^{\ell} + \check \beta_e^{\ell}) \nonumber  \\
&\quad + \sum_{e \in E} F_e \tilde x_e(1-\tilde d_e^s)(\delta_e^{\ell} + \eta_e^{\ell}) + \sum_{g \in G} \overline{P}_g  \tilde x_g(1-\tilde d_g^s)\zeta_g^{\ell},
\end{align}
where $L^s$ is the set of extreme points corresponding to the polyhedron characterized by dual constraints based on \eqref{psp} for primal variables $\mathbf f,\mathbf p, \mathbf q, \boldsymbol \theta$.

The constraint   $z(\tilde {\mathbf x}, \tilde{\mathbf d}^s) \le \varepsilon_{|s|}D$ should be satisfied for all $s \in S$.  Thus contingency feasibility conditions can be defined as follows,
\begin{align}\label{dsp3}
& \sum_{i \in V} D_i (\alpha_i^{s\ell} + \lambda_i^{\ell}) + \sum_{e \in E} M_e(1- x_e + \tilde d_e^s)(\hat \beta_e^{\ell} + \check \beta_e^{s\ell})  \nonumber\\
&+ \sum_{e \in E} F_e  x_e(1-\tilde d_e^s)(\delta_e^{\ell} + \eta_e^{\ell}) + \sum_{g \in G} \overline{P}_g   x_g(1-\tilde d_g^s)\zeta_g^{\ell} \\
 &\hskip 4.5cm \leq \varepsilon_{|s|} D, \  \forall \ell \in L^s, s \in S \nonumber
\end{align}

Since the objective is to minimize total investment cost and operating cost in the no-contingency state ($s=0$), we enforce all constraints for state 0 explicitly. For all other contingency states $s \in S\setminus \{0\}$, constraint set \eqref{dsp3} ensures that at least $\varepsilon_{|s|}$ fraction of the total demand is satisfied. The reformulation of \eqref{ef} is given as:

\begin{subequations}\label{rmp}
\begin{align}
\min_{\mathbf x, \mathbf f, \mathbf p, \mathbf q, \boldsymbol \theta}  &\  \sum_{e \in E} C_e x_e + \sum_{g\in G} C_gx_g + \sigma  \sum_{g\in G}  C_g^p p_g^0  \\
\textrm{s.t.} \ & \sum_{i \in V} D_i (\alpha_i^{\ell} + \lambda_i^{\ell}) + \sum_{e \in E} M_e(1- x_e + \tilde d_e^s)(\hat \beta_e^{\ell} + \check \beta_e^{\ell})  \nonumber\\
	&+ \sum_{e \in E} F_e x_e(1-\tilde d_e^s)(\delta_e^{\ell} + \eta_e^{\ell}) + \sum_{g \in G} \overline{P}_g  x_g(1-\tilde d_g^s)\zeta_g^{\ell} \nonumber \\
	&  \leq  \varepsilon_{|s|} D, \ \forall  \ell \in L^s, s \in S\setminus \{0\}  \label{rmp_f_cut2} \\
	&\sum_{g\in  G_i} p_g^0+ \sum_{e \in E_{.i}} f_{e}^{0} - \sum_{e \in E_{i.}} f_{e}^{0}  =  D_i, \ \forall i \in V  \\
    &B_{e}\big(\theta_{i_e}^{0} - \theta_{j_e}^{0}\big) - f_{e}^{0} + M_e(1-x_e) \geq  0, \ \forall e\in E \\
    &B_{e}\big(\theta_{i_e}^{0} - \theta_{j_e}^{0}\big) - f_{e}^{0} - M_e(1-x_e) \leq  0, \ \forall e\in E \\
	&-F_e \leq f_{e}^{0} \leq  F_{e}, \ \forall e \in  E   \\
	& 0 \leq p_g^{0} \leq \overline{P}_g x_g, \ \forall g \in G  \\
	&x_g \in \{0,1\}, \ \forall g\in G \label{pf_x_g}\\
	&x_e \in \{0,1\}, \ \forall e\in E \label{pf_x_e}
\end{align}
\end{subequations}

The number of constraints in formulation \eqref{ef} grows exponentially with the problem size, so we solve it via Benders Decomposition (BD).  At a typical iteration of BD, we consider the \emph{restricted master problem} (RMP) \eqref{rmp}, which has the same objective as \eqref{ef} but involves only a small subset of the constraints in \eqref{rmp_f_cut2}.  We briefly outline BD below.  For a detailed treatment of BD please refer to \cite{Benders1962}.

Let $t$ be the iteration number and let the initial RMP be problem \eqref{rmp} without any \eqref{rmp_f_cut2} constraints.  Let $\tilde {\mathbf x}^t$ be a concatenation of the expansion variables in the $t^{\textmd{th}}$ iteration.

\begin{algorithm}[H]
\caption{\emph{Benders Decomposition} (BD)}
\begin{algorithmic}[1]
\State $t \gets 0$
\State solve RMP
\State \textbf{if} RMP is infeasible
\State \hskip 0.4cm EXIT, TGEP is infeasible
\State \textbf{else}
\State \hskip 0.4cm Let $\tilde {\mathbf x}^t$ be the optimal solution to RMP
\State \hskip 0.4cm \textbf{for} {$s \in S$}, solve DSP$(\tilde {\mathbf x}^{t}, \tilde{\mathbf d}^s)$, let $z^t$ be the objective value
\State \hskip 0.8cm \textbf{if}  {$z^t> \varepsilon_{|s|} D$}
\State \hskip 1.2cm Add feasibility cut $(\ref{dsp3})$ to RMP
\State \hskip 0.8cm \textbf{end if}
\State \hskip 0.4cm \textbf{end for}
\State \hskip 0.4cm \textbf{if} no feasibility cut(s) added in Step 9
\State \hskip 0.8cm $\tilde {\mathbf x}^t$ is optimal, EXIT
\State \hskip 0.4cm \textbf{else}
\State \hskip 0.8cm $t \gets t+1$, go to Step 2
\State  \hskip 0.4cm \textbf{end if}
\State \textbf{end if}
\end{algorithmic}
\end{algorithm}

By using a Benders reformulation, we are able to decompose the extremely large formulation \eqref{ef} into a master problem and multiple subproblems (one for each contingency state).  In theory, this enables us to solve larger instances, which would not be possible by a direct solution of EF.  However, the extremely large number of contingency states makes direct application of Benders ineffective for large power systems and/or a non-trivial contingency budget (i.e., $k>1$). In the next section, we develop a custom cutting plane algorithm that evaluates all possible contingency states implicitly using a bilevel separation oracle.

\subsection{Online Contingency Scenario Generation}
Sizes of most power systems in operation (typically  thousands of generating units and transmission elements) may preclude direct solution of \eqref{ef}. Even using a decomposition algorithm (e.g. BD) may not be feasible because each contingency state must be considered explicitly. Our goal is to instead use a separation oracle that implicitly evaluates all contingency states and either identifies a violated one (a contingency with $j$ failures (for all $j=1,\cdots,k$) that cannot be survived by the current power system design) or provides a certificate that no such contingency state exists. If such a contingency exists, we use this contingency to generate a violated Benders feasibility cut, as described in the previous section, for the RMP. If no such contingency exists, then the current capacity expansion $\tilde {\mathbf x}$ is optimal and we terminate the algorithm.


\subsubsection{Power System Inhibition Problem (PSIP)}
Given a capacity expansion decision $\tilde {\mathbf x}$, the \emph{Power System Inhibition Problem} (PSIP) can be used  to determine the worst-case loss-of-load under any contingency with $j$ failures, for all $j=1,\cdots,k$. In this bilevel program, the upper level decisions $\mathbf d$ correspond to binary contingency selection decisions and the lower level decisions $(\mathbf f, \mathbf p, \mathbf q,\boldsymbol \theta)$ correspond to recourse power flow, generation scheduling, and load shedding decisions relative to the given contingency, prescribed by $\mathbf d$.

Note that in the prior model $\tilde {\mathbf{d}}^s$ was an input parameter,  whereas in this formulation, we are now selecting the elements of the contingency, with $\mathbf d$ becoming a vector of decision variables. For clarity of exposition, the superscript $s$ corresponding to variables $\textbf f, \textbf p, \textbf q$ and $\boldsymbol \theta$ has been removed, as the contingency state is not pre-specified, but rather part of the decision making process within the PSIP$(\tilde {\mathbf x}, j)$. PSIP$(\tilde {\mathbf x}, j)$ is given as follows:
\begin{subequations}\label{psip}
\begin{align}
\max_{\mathbf d}  \  \min_{\mathbf f,\mathbf p,\mathbf q,\boldsymbol \theta} \ \  &  \sum_{i \in V}  q_i   \label{psip_obj}\\
\textrm{s.t.}  \quad & \sum_{e \in  E} d_e + \sum_{g \in G} d_g = j, \label{psip_budget}\\
	 (\alpha_i)  \quad & \sum_{g\in  G_i} p_g  + \sum_{e \in E_{.i}} f_{e} - \sum_{e \in E_{i.}} f_{e} +q_i = D_i,  \  \forall i \label{psip_flow_consv}\\
	(\hat \beta_e) \quad &-B_{e}\big(\theta_{i_e} - \theta_{j_e}\big) + f_{e}  \leq M_e(1- \tilde x_e +  d_e), \ \forall e \label{psip_pa_ub}\\
	(\check \beta_e) \quad &B_{e}\big(\theta_{i_e} - \theta_{j_e}\big) - f_{e} \leq M_e(1- \tilde x_e +  d_e), \ \forall e  \label{psip_pa_lb} \\
	(\delta_e) \quad &f_{e} \leq  F_{e} \tilde x_e(1- d_e),  \ \forall e    \label{psip_flow_ub}\\
	(\eta_e) \quad &-f_{e} \leq  F_{e} \tilde x_e(1- d_e),  \ \forall e    \label{psip_flow_lb}\\
	(\zeta_g) \quad & 0 \leq p_g \leq \overline{P}_g \tilde   x_g (1-d_g), \ \forall g  \label{psip_gen_ub}  \\
	(\lambda_i) \quad & 0 \leq q_i \leq D_i, \ \forall i \label{psip_lol_ub}
\end{align}
\end{subequations}

The objective \eqref{psip_obj} is to maximize the minimum loss-of-load. For a given contingency state defined by $\mathbf d$, the objective of the power system operator (the inner minimization problem) is to determine the DCOPF such that the loss-of-load is minimized. Constraint \eqref{psip_budget} is a budget constraint limiting the number of power system elements that can be in the contingency. Constraints \eqref{psip_flow_consv} are standard flow conservation constraints. Constraints \eqref{psip_pa_ub} and \eqref{psip_pa_lb} together enforce Kirchhoff's voltage law, for active transmission elements. Constraints \eqref{psip_flow_ub} and \eqref{psip_flow_lb} are constraints associated with the capacity of each transmission element. Constraints \eqref{psip_gen_ub} limit the maximum capacity of each generating unit. If a generating unit $g$ is NOT part of the contingency (i.e.,  $d_g = 0$), then the maximum capacity of the generating unit is enforced, if the unit was added ($x_g=1$).  Otherwise, the power output of the generating unit must be zero.

The upper-level decisions of this bilevel program are to select a contingency, using the binary variables $\mathbf d$, that maximizes the subsequent loss-of-load in the lower-level problem.


\subsubsection{A MILP Reformulation of PSIP}
Bilevel programs like \eqref{psip} cannot be solved directly. One approach is to reformulate the bilevel program by dualizing the inner minimization problem.  For fixed values of $\mathbf d$, the inner minimization problem is a linear program that is always feasible. By strong duality, and combining the upper level of \eqref{psip} we can reformulate the bilevel program as a single level bilinear program, where the objective function of the dualized problem contain terms associated with the product of the upper-level contingency selection variables $\mathbf d$ and the lower-level dual variables $(\hat {\boldsymbol  \beta}, \check {\boldsymbol  \beta},\boldsymbol \delta,\boldsymbol  \eta, \boldsymbol \zeta)$ associated with flow balance, transmission flow, transmission capacity, and generation capacity constraints.  However, with additional variables and constraints, these bilinear terms can be linearized.

Each bilinear term can be linearized using the following strategy.  Let $u\le 0$ and $v \le 0$ be continuous variables and $b \in \{0,1\}$.  Then the bilinear term, $bu$, can be linearized as follows. Letting $v = bu$, we introduce the following three constraints to linearize the bilinear term $bu$.
\begin{subequations}\label{cons_linearize}
\begin{align}
	& v \ge u-U(1-b) \label{lin1}\\
	& v \ge -Ub \label{lin2}\\
	& v \le u + U(1-b) \label{lin3}
\end{align}
\end{subequations}

Here, parameter $U$ represents a valid upper bound for continuous variable $u$ and satisfies $U \ge |u|$. Assessing these three constraints for both binary values of $b$ show that they indeed provide a valid linearization.  If $b=0$, then constraints  (\ref{lin2}) and $v\le 0$ together imply that $v=0$.  With $v=0$, constraints  (\ref{lin1}) and (\ref{lin3}) together imply that $-U \le u\le U$, which are never binding. If $b=1$, then constraints (\ref{lin1}) and (\ref{lin3}) together imply  $u=v$ and constraints (\ref{lin2}) and $v\le 0$   implies $-U \le v \le 0$, which is valid.

We follow a similar strategy to linearize all five bilinear terms  $(\hat {\boldsymbol  \beta} d_e, \check {\boldsymbol  \beta}d_e,\boldsymbol \delta d_e,\boldsymbol  \eta d_e, \boldsymbol \zeta d_g)$.  Define continuous variables $(\mathbf r^1, \mathbf r^2,\mathbf  r^3, \mathbf r^4, \mathbf r^5)$ and let $r_e^1 =  \hat \beta_e d_e$, $r_e^2 = \check \beta_e d_e $, $r_e^3 =  \delta_e d_e$, $r_e^4 = \eta d_e$, and $r_g^5 = \zeta_g d_g$.  Following the same linearization strategy introduced in (\ref{cons_linearize}), we now state the full mixed integer linear PSIP formulation for completeness, which we call the  {\em Mixed-Integer Power System Inhibition Problem}, M-PSIP$(\tilde {\mathbf x}, j)$:

\begin{subequations}\label{psip2}
\begin{align}
\hskip -0.2cm \max_{\mathbf d, \boldsymbol \alpha, \hat {\mathbf \beta}, \check {\mathbf \beta}, \delta, \eta, \zeta, \mathbf r} \ & \sum_{i\in V} D_i (\alpha_i + \lambda_i) + \sum_{e \in E} M_e(1- \tilde x_e)(\hat \beta_e + \check \beta_e) \nonumber \\
&\hskip -1.4cm + \sum_{e \in E} \Big (M_e(r_e^1 + r_e^2) +  F_e   \tilde  x_e(\delta_e + \eta_e) - F_e   \tilde  x_e(r_e^3 + r_e^4)  \Big)\nonumber \\
&\hskip -1.4cm + \sum_{g \in G} \Big( \overline{P}_g   \tilde  x_g\zeta_g  -  \overline{P}_g   \tilde   x_g r_g^5 \Big) \\
&\hskip -1.8cm \textrm{s.t.}\quad \sum_{e \in  E} d_e + \sum_{g \in G} d_g = j, \label{psip2_budget}\\
 &\hskip -1cm  \alpha_{j_e}-\alpha_{i_e}+\hat \beta_e - \check \beta_e +\delta_e-\eta_e=0, \ \forall e\\
  &\hskip -1cm\alpha_{i_g}+\zeta_g\leq 0, \ \forall g\\
 &\hskip -1cm\alpha_i + \lambda_i \leq 1, \ \forall i\\
 &\hskip -1cm\sum_{e\in E_{i_e=i,.}} B_e (\check \beta_e - \hat \beta_e) + \sum_{e\in E_{.,j_e=i}} B_e (\hat \beta_e - \check \beta_e) =0, \ \forall i  \\
&\hskip -1cm\hat \beta_e \leq 0, \ \check \beta_e \leq 0, \ \delta_e\le 0, \ \eta_e\le 0, \ \forall e  \\
&\hskip -1cm r_e^1 \leq 0, \  r_e^2 \leq 0, \  r_e^3\le 0, \  r_e^4\le 0, \ \forall e  \\
&\hskip -1cm\zeta_g\le 0,\ r_e^5 \le 0, \ \forall g \\
&\hskip -1cm \lambda_i \le 0, \ \forall i
\end{align}
\end{subequations}


Next, we outline an algorithm for \emph{optimally} solving problem \eqref{ef} that combines a Benders decomposition with the aid of an oracle given by \eqref{psip2}, which acts as a separation problem.  A given capacity expansion  $\tilde {\mathbf x}$ is optimal if the oracle cannot find a contingency of size $j$, for any $j=1,\cdots,k$, that results in a loss-of-load above the allowable threshold $\varepsilon_j D$.

For each contingency budget $j$, we can check for $j$-element contingencies by solving M-PSIP using a failure budget of $j$ (i.e. the right-hand side of inequality (\ref{psip2_budget}) is set to $j$).  Whenever the oracle determines that the capacity expansion decision $\tilde {\mathbf x}$ is \emph{not} $N$-$k$-$\boldsymbol \varepsilon$ compliant,  it returns a contingency $\mathbf d$ that results in a loss-of-load, above the allowable threshold $\varepsilon_j D$ for $j$-element failures.

Let $r$ be the iteration number and let the initial RMP be problem \eqref{rmp} without any  \eqref{rmp_f_cut2} constraints.  Let $\tilde {\mathbf x}^r$ be a concatenation of the expansion variables $(x_e^r,x_g^r)$.

\begin{algorithm}[H]
\caption{\emph{Online Contingency Screening} (OCS)}
\begin{algorithmic}[1]
\State $t \gets 0$
\State Solve RMP for iteration $t$
\State \textbf{if} RMP is infeasible
\State \hskip 0.4cm EXIT, TGEP is infeasible
\State \textbf{else}
\State \hskip 0.4cm Let $\tilde {\mathbf x}^t$ be the optimal solution to RMP
\State  \hskip 0.4cm \textbf{for} $j=1,\cdots, k$
\State \hskip 0.8cm Solve M-PSIP$(\tilde {\mathbf x}^t, j)$, let $z^t_j$ be the objective value\\
\hskip 0.8cm  and  $\mathbf d_j^t$ be the contingency selection decision
\State \hskip 0.8cm \textbf{if} { $z_j^t > \varepsilon_j D$} then 
\State \hskip 1.2cm Solve DSP$(\tilde {\mathbf x}^t, \mathbf d_j^t)$, add feasibility cut \eqref{dsp3} to RMP
\State \hskip 1.2cm $t \gets t+1$, go to step 2
\State \hskip 0.8cm  \textbf{end if}
\State \hskip 0.4cm \textbf{end for}
\State \textbf{end if}
\State $\tilde {\mathbf x}^t$ is optimal, EXIT
\end{algorithmic}
\end{algorithm}

At each iteration, either a contingency that results in loss-of-load above the allowable threshold is identified and a corresponding feasibility cut is generated and added to RMP, or no contingencies are found, which means that the current solution is optimal and the algorithm can terminate.


\section{Numerical Experiments}\label{sec4}
We implemented the proposed models and algorithms in C++ and CPLEX 12.1 via ILOG Concert Technology 2.9. All experiments were run on a machine with four quad-core 2.93G Xeon with 96G of memory.  For the following computational experiments, a single CPU and up to 8GB of RAM were allocated. The optimality gap was set to be 0.1\% for CPLEX.

We have tested our models and algorithms on the IEEE 30-bus and IEEE 57-bus systems \cite{TestData}. For each power system, we consider five different contingency budgets $k=0,1,2,3,$ and $4$.
Altogether, we consider 10 instances.


Table \ref{tab1} compares the run times for the three different approaches. For each of the 10 instances, $m$ provides the number of distinct contingencies. Initially for each test system, we replicate a subset of existing generating units and transmission lines to create a set of candidate elements. With these candidate elements as a starting point, we iteratively solve the PSIP problem for $k$ and $\varepsilon$ values presented in Table~\ref{tab1} using OCS.  In all 10 instances, the values of $\varepsilon_i$  are invariant. That is, for instance, $\varepsilon_2=0.05$   when $k$ was chosen to be 2,3, or 4. For these, we identify vulnerabilities in the power system and introduce additional candidate generation and transmission elements. We follow this method to create the augmented the IEEE 30-bus and IEEE 57-bus test systems for the computational experiments presented subsequently. We want to note that our contributions in this paper are in algorithmic fundamentals, and  thus are not very sensitive to
 the particular problem instances. Hence, the particular  list of candidates lines  for a problem instance will only have a minor effect on the performances of our algorithms.


 \begin{table}[th]
\caption{Run times for different solution approaches}
\centering
\begin{tabular}{c c c c| c c c}
\hline \hline
\multicolumn{4}{c|}{} & \multicolumn{3}{c}{Solution time (secs)} \\[0.5ex]
Test Systems 		&$m$ 		& $k$	&  $\boldsymbol \varepsilon$ 		& EF        		& BD        	 & OCS           	 \\[0.5ex]
\hline
IEEE 30-bus		&	0		&	0	&	0 		    				& 0 			& 0 			&0			 \\
         				&	152		&	1	&	0 		    				& 76 			& 2 			& 4				    \\
            			&	$>11K$	&	2	&	0.05 						& x 			& 87 			& 16				 \\
        				&	$>500K$	&	3	&	0.10 						& x 			& 3,126 		& 128				 \\
        				&	$>21M$	&	4	&	0.20 						& x 			& x 			& 141				 \\ [0.5ex]
\hline
IEEE 57-bus     		&	0		&	0	&	0 		    				& 0 			& 0 			& 0				    \\
         				&	110		&	1	&	0 		    				& 27 			& 164 		&120				    \\
            			&	$>5K$	&	2	&	0.05 						& x 			& 870 		& 36				 \\
        				&	$>200K$	&	3	&	0.10 						& x 			& x 			& 51				 \\
        				&	$>5M$	&	4	&	0.20 						& x 			& x 			& 65				 \\ [0.5ex]
\hline\hline
\end{tabular}
\label{tab1}
\end{table}

Table \ref{tab1} provides the run time (in CPU seconds) for each instance under the three different approaches.  In this table ``x'' means the algorithm failed to complete at the end of 2 hours. Note that the first approach, the extensive form (EF), can only solve the smallest of instances. This is because of the sheer size of the problem, in which, for each contingency, a full DCOPF problem must be embedded in the formulation. As the number of contingencies grows, this formulation quickly becomes intractable. Note that  we are working with small data sets here and target systems will be in the order of thousands  of elements,  which will make EF intractable even sooner.

The second approach, BD, bypasses this problem via a Benders decomposition, with corresponding delayed cut generation. However, this still suffers from the combinatorial growth in the number of contingency states -- for each contingency, a subproblem (DSP) must be solved to check for violated feasibility cuts to add to the RMP. We see that larger problem instances can be solved, relative to EF, but the BD approach nonetheless cannot solve the largest problem instances.

With the OCS approach, we see that all instances of the problem can be solved, in all cases in under three minutes and frequently in only a few seconds. This is a result of the combination of the strength of the Benders cuts, enabling the problem to be solved in a very limited number of iterations, and also the fact that we are able to implicitly evaluate the contingencies in order to identify a violated contingency and then quickly find its corresponding feasibility cut by solving a single linear program (DSP).

Table \ref{tab2} provides us with further evidence  for the scalability of OCS. For each instance, we see the total number of possible contingency states $m$ and then the number of contingency states for which corresponding feasibility cuts were actually generated, denoted by `cont.' in the table. Clearly, it is a very tiny fraction of the possible number of contingencies, which is critical to the tractability of the approach. The remaining columns of this table breakdown the total run time by time spent on the three components of the algorithm -- the RMP, which identifies a candidate network design; the mixed-integer linear power system inhibition problem (M-PSIP), which identifies a contingency that cannot be overcome by the current network design; and the dual subproblems (DSP), which generates the feasibility cuts.

\begin{table}[th]
\caption{OCS runtime breakdown}
\centering
\begin{tabular}{c c c c c c c c}
\hline \hline
Test Systems 	&$m$ 		& $k$	&  $\varepsilon$ 	&RMP        	&M-PSIP       	&DSP     &cont.   	\\[0.5ex]
\hline
IEEE 30-bus	&	0		&	0	&	0 		    	& 0 		    	& 0		        & 0			&0	 \\
         			&	152		&	1	&	0 		    	& 1 			& 3 			& 0			&2	 \\
        			&	$>11K$	&	2	&	0.05 			& 1 			& 15 			& 0			&5		 \\
        			&	$>500K$	&	3	&	0.10 			& 1 			& 127 			& 0		&8		 \\
        			&	$>21M$	&	4	&	0.20 			& 1 			& 140		    & 0		&10			 \\ [0.5ex]
\hline
IEEE 57-bus     &	0		&	0	&	0 		    	& 0 		    	& 0		        & 0			&0	 \\
       			&	110		&	1	&	0 		    	& 83			& 37 			& 0			&3		 \\
        			&	$>5K$	&	2	&	0.05 			& 15 			& 21 			& 0			&4		 \\
        			&	$>200K$	&	3	&	0.10 			& 14 			& 37 			& 0			&7	 \\
        			&	$>5M$	&	4	&	0.20 			& 14 			& 51		    & 0			&8 \\ [0.5ex]
\hline \hline
\end{tabular}
\label{tab2}
\end{table}


\section{Conclusion}\label{sec5}
We studied  the transmission and generation expansion problem with contingency constraints.  More specifically,  we investigated the problem of improving an electric power system at a minimum cost   by adding generators and transmission lines, such that it satisfies the $N$-$k$-$\varepsilon$ survivability criterion. This survivability criterion is a generalization of the well-known $N$-$k$  criterion, and it requires that at least $(1-\varepsilon_j)$  fraction of the total demand is met even after failures of any  $j$  system components, for all $j=1,\cdots,k$. This design problem  adds another level of complexity to the contingency analysis problem, since the contingency analysis is only one of the constraints  in  the design optimization problem.   We proposed two algorithms: one is based on the Benders decomposition approach, and the other is based on online contingency state generation. The latter approach avoids the combinatorial explosion by seeking vulnerabilities in the current solution, and generating constraints to exclude such infeasible solutions.   We tested our proposed approaches on the IEEE 30-bus and the IEEE 57-bus systems.  Computational results show the proposed online contingency generation algorithm, which uses a bilevel separation technique to implicitly consider all exponential number of contingencies, significantly outperforms a standard Benders decomposition.  We were able to solve all instances in our experiment in under three minutes, while the extensive form and the Benders Decomposition algorithm  failed to complete at the end of 2 hours.

We believe that  this paper will provide the fundamentals for many other studies in contingency-aware transmission and generation expansion.  As an example, we want to apply for methods to full-scale systems. While
our results are very promising in terms of scalability, full-scale problems will surely pose some computational challenges  and will require adopting high-performance computing resources.  Also our current model assumes
all  failures happen simultaneously. In order to reflect practical operation situations, where failures may happen consecutively, new models that consider timing between system element failures are needed. Additionally, unit commitment and de-commitment is not considered in our current model.  We plan to extend our models for these cases. Finally, we worked with a deterministic model,  and it is essential, to take stochasticity into account for planning problems.  We believe our current frame work can be naturally extended  for stochastic problems.

\section* {Acknowledgements.} This work was funded by the applied mathematics program at the United States Department of Energy and by the Laboratory Directed Research \& Development (LDRD) program at Sandia National Laboratories. Sandia National Laboratories is a multiprogram laboratory operated by Sandia Corporation, a wholly owned subsidiary of Lockheed Martin Corporation, for the United States Department of Energy's National Nuclear Security Administration under contract DE-AC04-94AL85000.

\ifCLASSOPTIONcaptionsoff
  \newpage
\fi

\bibliographystyle{IEEE}
\bibliography{inhibition}

\begin{thebibliography}{10}

\bibitem{NERC}
North American Electric Reliability~Corporation (NERC),
\newblock ``Standard tpl-001-1--system performance under normal conditions'',
  Feb. 2011.

\bibitem{Donde08}
V.~Donde, V.~Lopez, B.~Lesieutre, A.~P{\i}nar, C.~Yang, and J.~Meza,
\newblock ``\href{http://csmr.ca.sandia.gov/~apinar/papers/IEEE-TPS.pdf}{Severe
  Multiple Contingency Screening in Electric Power Systems}'',
\newblock {\em IEEE Transactions on Power Systems}, vol. 23, pp. 406--417,
  2008.

\bibitem{lesieutre08}
B.C. Lesieutre, A.~Pinar, and S.~Roy,
\newblock ``Power system extreme event detection: The vulnerability frontier'',
\newblock in {\em Proc. 41st Hawaii International Conference on System
  Sciences}, Hawaii, 2008.

\bibitem{Pinar2010}
A.~Pinar, J.~Meza, V.~Donde, and B.~Lesieutre,
\newblock ``Optimization strategies for the vulnerability analysis of the
  electric power grid'',
\newblock {\em SIAM J. Optim.}, vol. 20, pp. 1786--1810, 2010.

\bibitem{Bienstock2010}
D.~Bienstock and A.~Verma,
\newblock ``The n-k problem in power grids: new models, formulations, and
  numerical experiments'',
\newblock {\em SIAM J. Optim.}, vol. 20, pp. 2352--2380, 2010.

\bibitem{Arroyo2010}
J.M. Arroyo,
\newblock ``Bilevel programming applied to power system vulnerability analysis
  under multiple contingencies'',
\newblock {\em IET Gener Transm. Distrib.}, vol. 4, pp. 178--190, 2010.

\bibitem{SalWB04}
J.~Salmeron, K.~Wood, and R.~Baldick,
\newblock ``Analysis of electric grid security under terrorist threat'',
\newblock {\em IEEE Trans. Power Syst.}, vol. 19, pp. 905--912, 2004.

\bibitem{SalWB09}
J.~Salmeron, K.~Wood, and R.~Baldick,
\newblock ``Worst-case interdiction analysis of large-scale electric power
  grids'',
\newblock {\em IEEE Trans. Power Syst.}, vol. 24, pp. 96--104, 2009.

\bibitem{Fan2011}
N.~Fan, H.~Xu, F.~Pan, and P.M. Pardalos,
\newblock ``Economic analysis of the n-k power grid contingency selection and
  evaluation by graph algorithms and interdiction methods'',
\newblock {\em Energy Syst.}, vol. 2, pp. 313--324, 2011.

\bibitem{Street2011}
A.~Street, F.~Oliveira, and J.M. Arroyo,
\newblock ``Contingency-constrained unit commitment with $n$-$k$ security
  criterion: A robust optimization approach'',
\newblock {\em IEEE Trans. Power Syst.}, vol. 26, pp. 1581--1590, 2011.

\bibitem{CheM2005}
Q.~Chen and J.D. McCalley,
\newblock ``Identifying high risk n-k contingencies for online security
  assessment'',
\newblock {\em IEEE Trans. Power Syst.}, vol. 20, pp. 823--834, 2005.

\bibitem{CarL2002}
B.~A. Carreras, V.~E. Lynch, I.~Dobso, and D.~E. Newma,
\newblock ``Critical points and transitions in an electric power transmission
  model for cascading failure blackouts'',
\newblock {\em Chaos}, vol. 12, pp. 985--994, 2002.

\bibitem{LiC2009}
F.~Li and A.~Chegu,
\newblock ``High order contingency selection using particle swarm optimization
  and tabu search'',
\newblock in {\em Proc. of the 15th Inter. Conf. on Intelligent System
  Applications to Power Systems}, Curitiba, Brazil, Nov. 2009.

\bibitem{HinC2010}
P.~Hines, E.~Cotilla-Sanchez, and S.~Blumsack,
\newblock ``Do topological models provide good information about electricity
  infrastructure vulnerability?'',
\newblock {\em Chaos}, vol. 20, pp. 033122, 2010.

\bibitem{MorG2001}
H.~Mori and Y.~Goto,
\newblock ``A tabu search based approach to (n-k) static contingency selection
  in power systems'',
\newblock {\em Proc. 2001 IEEE Inter. Conf. on Systems, Man, and Cybernetics},
  vol. 3, pp. 1954--1959, 2001.

\bibitem{CheJ2009}
Y.~Chen, S.~Jin, D.~Chavarria-Miranda, and Z.~Huang,
\newblock ``Application of cray xmt for power grid contingency selection'',
\newblock in {\em Proc. Cray User Group 2009}, Atlanta, GA, May 2009.

\bibitem{Latorre2003}
G.~Latorre, R.D. Cruz, J.M. Areiza, and A.~Villegas,
\newblock ``Classification of publications and models on transmission expansion
  planning'',
\newblock {\em IEEE Trans. Power Syst.}, vol. 18, pp. 938--946, 2003.

\bibitem{nsw11}
R.~Chen, A.~Cohn, and A.~Pinar,
\newblock ``\href{http://arxiv.org/abs/1109.1801}{An Implicit Optimization
  Approach for Survivable Network Design}'',
\newblock in {\em Proc. 2011 IEEE $1^{st}$ International Network Science
  Workshop (NSW 2011)}, 2011.

\bibitem{Pmaps11}
R.~Chen, A.~Cohn, N.~Fan, and A.~Pinar,
\newblock ``$n$-$k$-$\epsilon$ survivable power system design'',
\newblock in {\em $12^{th}$ International Conference on Probabilistic Methods
  Applied to Power Systems (PMAPS12)}, Istanbul, Turkey, June 2012.

\bibitem{Delgadillo2011}
A.~Delgadillo, J.M. Arroyo, and N.~Alguacil,
\newblock ``Power syustem defense planning against multiple contingencies'',
\newblock in {\em Proc. 17th Power Syst. Compt. Conf. (PSCC)}, Stochholm,
  Sweden, Aug. 2011.

\bibitem{Romero2012}
N.~Romero, N.~Xu, L.K. Nozick, I.~Dobson, and D.~Jones,
\newblock ``Investment planning for electric power systems under terrorist
  threat'',
\newblock {\em IEEE Trans. Power Syst.}, vol. 27, pp. 108--116, 2012.

\bibitem{Carrion2007}
M.~Carrion, J.M. Arroyo, and N.~Alguacil,
\newblock ``Vulnerability-constrained transmission expansion planning: A
  stochastic programming approach'',
\newblock {\em IEEE Trans. Power Syst.}, vol. 22, pp. 1436--1445, 2007.

\bibitem{Choi2007}
J.~Choi, T.D. Mount, and R.J. Thomas,
\newblock ``Transmission expansion planning using contingency criteria'',
\newblock {\em IEEE Trans. Power Syst.}, vol. 22, pp. 2249--2261, 2007.

\bibitem{Moulin2010}
L.~Moulin, M.~Poss, and C.~Sagastiz\'{a}bal,
\newblock ``Transmission expansion planning with re-design'',
\newblock {\em Energy Syst.}, vol. 1, pp. 113--139, 2010.

\bibitem{Zhang2011}
H.~Zhang, V.~Vittal, G.T. Heydt, and J.~Quintero,
\newblock ``A mixed-integer linear programming approach for multi-stage
  security-constrained transmission expansion planning'',
\newblock {\em IEEE Trans. Power Syst.}, 2011.

\bibitem{Jin2011}
S.~Jin, S.M. Ryan, J.-P. Watson, and D.L. Woodruff,
\newblock ``Modeling and solving a large-scale generation expansion planning
  problem under uncertainty'',
\newblock {\em Energy Syst.}, vol. 2, pp. 209--242, 2011.

\bibitem{Bent2011}
R.~Bent, A.~Berscheid, and A.L. Toole,
\newblock ``Generation and transmission expansion planning for renewable energy
  integration'',
\newblock in {\em Proc. of Power Syst. Comp. Conf. (PSCC)}, Stockholm, Sweden,
  Aug. 2011.

\bibitem{Weijde2010}
A.H. van~der Weijde and B.F. Hobbs,
\newblock ``Transmission planning under uncertainty: A two-stage stochastic
  modelling approach'',
\newblock in {\em Proc. 7th Int. Conf. Euro. Energy Market (EEM)}, Madrid,
  Spain, 2010.

\bibitem{Benders1962}
J.F. Benders,
\newblock ``Partitioning procedures for solving mixed-variables programming
  problems'',
\newblock {\em Numerische Mathematik}, vol. 10, pp. 237--260, 1962.

\bibitem{TestData}
IEEE reliability~test data,
\newblock '',
\newblock http://www.ee.washington.edu/research/pstca/.

\end{thebibliography}

\end{document}